\documentclass[12pt]{article}
\usepackage[intlimits]{amsmath}
\usepackage{amssymb}
\usepackage{mathrsfs}
\usepackage[dvips]{graphicx}
\usepackage{longtable}

\title{Counterexamples to Borsuk's conjecture on spheres of small radii\footnote{This work is done 
under the financial support of the following grants: the grant 09-01-00294 of 
Russian Foundation for Basic Research,
the grant MD-8390.2010.1 of the Russian President, the grant NSh-8784.2010.1 supporting Leading scientific schools of Russia.}}

\author{A.B. Kupavskii\footnote{Moscow State University, Mechanics and Mathematics Faculty, Department of Number Theory.},
A.M. Raigorodskii\footnote{Moscow State University, Mechanics and Mathematics Faculty, Department of Mathematical Statistics 
and Random Processes; Moscow Institute of Physics and Technology, Faculty of Innovations and High Technology, Department of Data 
Analysis; Yandex research laboratories.}}

\date{}

\DeclareMathOperator{\diam}{diam}

\begin{document}

\maketitle

%\begin{large}

\begin{abstract}

In this work, the classical Borsuk conjecture is discussed, which states that any set of diameter 1 in the 
Euclidean space $ {\mathbb R}^d $ can be divided into $ d+1 $ parts of smaller diameter. During the last two decades, 
many counterexamples to the conjecture have been proposed in high dimensions. However, all of them are sets of diameter 1 that lie 
on spheres whose radii are close to the value $ \frac{1}{\sqrt{2}} $. The main result of this paper is as follows: {\it for any $ r > \frac{1}{2} $, 
there exists a $ d_0 $ such that for all $ d \ge d_0 $, a counterexample to Borsuk's conjecture can be found on a sphere $ S_r^{d-1} \subset 
{\mathbb R}^d $.}  

\end{abstract}

\section{Introduction}

This paper is devoted to the classical Borsuk partition problem. In 1933 K. Borsuk (see \cite{Bor1}) posed the following question: {\it is it true
that any set $ \Omega \subset {\mathbb R}^d $ having diameter 1 can be divided into some parts $ \Omega_1, \dots, \Omega_{d+1} $ whose diameters
are strictly smaller than 1?} Here by the {\it diameter} of a set $ \Omega $ we mean the quantity
$$
\diam \Omega = \sup_{{\bf x},{\bf y} \in \Omega} |{\bf x}-{\bf y}|, 
$$
where, in turn, $ |{\bf x}-{\bf y}| $ denotes the standard Euclidean distance between vectors. 

Let us use some additional notation. So by $ f(\Omega) $ we denote the value
$$
f(\Omega) = \min \{f: ~ \Omega = \Omega_1 \cup \ldots \cup \Omega_f, ~~ \forall~ i~~ \diam \Omega_i < \diam \Omega\}.
$$
Furthermore, 
$$
f(d) = \max_{\Omega \subset {\mathbb R}^d,~ \diam \Omega = 1} f(\Omega), 
$$
i.e., $ f(d) $ is the minimum number of parts of smaller diameter, into which an arbitrary set of diameter 1 in $ {\mathbb R}^d $ can be divided. In these 
terms, Borsuk's question is as follows: {\it is it true that always $ f(d) = d+1 $?} The positive answer on this question is usually 
called ``Borsuk's conjecture''. 

The history of Borsuk's conjecture is somewhat dramatic. Almost all the specialists in the field of combinatorial geometry strongly believed that 
the conjecture should be true. Multiple results supporting it have been proved. For example, if $ \Omega $ has smooth boundary, then one 
certainly obtains  
$ f(\Omega) \le d+1 $. This result is due to H. Hadwiger (see \cite{Had}), and it seems to show the evidence of the conjecture. However, in 1993 
J. Kahn and G. Kalai published a break-through paper, where they constructed a {\it finite} set of points in a very high dimension $ d $ that could not 
be decomposed into $ d+1 $ subsets of smaller diameter (see \cite{KK}). 

Now, Borsuk's conjecture is shown to be true for $ d \le 3 $ and false for $ d \ge 298 $ (see \cite{Rai1}). Also, we know that 
$$
(1.2255...+o(1))^{\sqrt{d}} \le f(d) \le (1.224...+o(1))^d.
$$
Here the lower bound was found in 1999 by A.M. Raigorodskii (see \cite{Rai2}) and the upper estimate is due to O. Schramm (see \cite{Sch}). 

The colorful history of Borsuk's conjecture is exhibited in numerous books and survey papers. We refer the reader to \cite{Rai1}, \cite{Gr}, 
\cite{BMS}, \cite{BMP}, \cite{Rai3}, \cite{Rai4}.    

\section{Formulation of the problem and statements of the results}

A careful analysis of all the known counterexamples to Borsuk's conjecture shows that they are always finite sets of points in 
$ {\mathbb R}^d $ lying on spheres whose radii are close to $ \frac{1}{\sqrt{2}} $. This is 
quite natural, since, by Jung's theorem (see \cite{Jung}), 
any set in $ {\mathbb R}^d $ having diameter 1 can be covered by a ball of radius $ \sqrt{\frac{d}{2d+2}} \sim 
\frac{1}{\sqrt{2}} $, and the intuition is that in order to get a counterexample, we have to take a set with as big covering ball as possible.   
The main result of our work is the following theorem which completely breaks such intuition. 

\vskip+0.2cm

\noindent {\bf Theorem 1.} {\it Let $ S_r^{d-1} \subset {\mathbb R}^d $ be the sphere of radius $ r $ with centre at the origin. 
For any $ r > \frac{1}{2} $, there exists a $ d_0 = d_0(r) $ such that for every $ d \ge d_0 $, 
one can find a set $ \Omega \subset S_r^{d-1} $ which has diameter 1 and does not admit a partition into $ d+1 $ parts of smaller diameter.} 

\vskip+0.2cm

Actually, we shall prove an even stronger result. In order to formulate it, let us introduce the quantity 
$$
f_r(d) = \max_{\Omega \subset S_r^{d-1},~ \diam \Omega = 1} f(\Omega). 
$$
In these terms, Theorem 1 sais that for any $ r > \frac{1}{2} $, there exists a $ d_0 = d_0(r) $ such that for every $ d \ge d_0 $, 
$ f_r(d) > d+1 $. Moreover, one has 

\vskip+0.2cm

\noindent {\bf Theorem 2.} {\it For any $ r > \frac{1}{2} $, there exist numbers $ k = k(r) \in {\mathbb N} $, $ c = c(r) > 1 $ and 
a function $ \delta = \delta(d) = o(1) $ such that 
$$ 
f_r(d) \ge (c+\delta)^{\sqrt[2k]{d}}.
$$}

\vskip+0.2cm

Theorem 2 means that if $ r $ is fixed and exceeds $ \frac{1}{2} $, then the order of magnitude of the value $ f_r(d) $ is at least 
$ e^{g(d)} $, where $ g(d) $ is just a fixed positive power of $ d $. So not only the quantity $ f_r(d) $ is greater than $ d+1 $
starting from some $ d_0 $, but also it is {\it substantially} greater than the conjectured value. This fact allows us to use some 
optimization and to prove eventually the following theorem. 

\vskip+0.2cm

\noindent {\bf Theorem 3.} {\it Let $ r = r(d) = \frac{1}{2} + \varphi(d) $, where $ \varphi = o(1) $ and $ \varphi(d) \ge 
c \frac{\ln\ln d}{\ln d} $ for all $ d $ and a large enough $ c > 0 $. Then, there exists a $ d_0 $ such that for $ d \ge d_0 $, 
$ f_{r(d)}(d) > d+1 $.}

\vskip+0.2cm

In other words, we can disprove Borsuk's conjecture by constructing sets of diameter 1 that lie 
on spheres with radii tending to $ \frac{1}{2} $ with $ d \to \infty $. Here, in the case of $ \varphi(d) = \Theta\left(\frac{\ln\ln d}{\ln d}\right) $, 
there is already no room to spare in the bound $ f_r(d) > d+1 $. So using our approach, one only may discuss the value of a constant $ c $ in 
Theorem 3, which is of course not significant. 

Now, a natural question arises: perhaps $ f_r(d) \le d+1 $, provided $ r = r(d) = O\left(\frac{\ln\ln d}{\ln d}\right) $? Unfortunately, 
we only can prove 

\vskip+0.2cm

\noindent {\bf Theorem 4.} {\it Let $ r = r(d) = \frac{1}{2} + \varphi(d) $, where $ \varphi = O(1/d) $. Then, $ f_r(d) \le d+1 $.} 

\vskip+0.2cm

Thus, we get a gap between the two functions $ \frac{1}{d} $ and $ \frac{\ln\ln d}{\ln d} $. To reduce this gap would be of a great 
interest. 

\vskip+0.2cm

The structure of the rest of our paper will be as follows. In Section 3, we shall prove both Theorem 1 and Theorem 2. Section 4 will be 
devoted to the proof of Theorem 3. In Section 5, we shall briefly discuss Theorem 4 which is in fact very simple. In Section 6, we shall show that 
in some sense, the constructions from Sections 3 and 4 are best possible. In Section 7, we shall say a few words about possible extensions 
and strengthenings of our main results.  

\section{Proofs of Theorems 1 and 2}

Since Theorem 1 is an immediate consequence of Theorem 2, we just prove the second result. 

\vskip+0.1cm

\noindent Let us fix an arbitrary $ r > \frac{1}{2} $. Without loss of generality, we also assume that $ r < \frac{1}{\sqrt{2}} $. 
Let $ d $ be large enough (during the proof, we will see what it means). To make clear our further exposition, we subdivide it into 
four parts. So in Subsection 3.1, we construct a $ d' $-dimensional set $ \Omega' $ with $ d' < d $; in Subsection 3.2, we show that $ \Omega' $ can 
be transformed into some $ \Omega \subset S_r^{d-1} $; in Subsection 3.3, we prove that $ f_r(d) \ge 
f(\Omega) \ge f(\Omega') \ge (c+\delta)^{\sqrt[2k]{d}} $
with some appropriate $ c > 1 $ and $ \delta = o(1) $; in Subsection 3.4, we prove a key lemma which is formulated in Subsection 3.3.

\subsection{Construction}
  
Take 
$$
k = \min \left\{k' \in {\mathbb N}: ~ r^2 > \frac{2k'+1}{8k'}\right\}. 
$$
This value is correctly defined, since $ r > \frac{1}{2} $ and $ \frac{2k'+1}{8k'} \to \frac{1}{4} $ as $ k' \to \infty $.   
Put
$$
n = \max\left\{m: ~ m \equiv 0 \pmod 4, ~~ m^{2k} < d\right\}. 
$$
Clearly 
$$ 
\sqrt[2k]{d} - 5 \le n < \sqrt[2k]{d}. 
\eqno{(1)}
$$

Consider the function
$$
u(a_0) = \frac{1+2k\left(\frac{a_0}{2}\right)^{2k-1}}{2+4k \left(\frac{a_0}{2}\right)^{2k-1} + 
(4k-2) \left(\frac{a_0}{2}\right)^{2k}}.
\eqno{(2)}
$$
Obviously, 
$$ 
u(0) = \frac{1}{2}, ~~~ u(2) = \frac{2k+1}{8k} < r^2 < \frac{1}{2}, 
$$
and $ u $ is monotone decreasing on the interval $ (0,2) $. Therefore, the equation 
$$
\frac{1+2k\left(\frac{a_0}{2}\right)^{2k-1}}{2+4k \left(\frac{a_0}{2}\right)^{2k-1} + 
(4k-2) \left(\frac{a_0}{2}\right)^{2k}} = r^2
$$
has a unique solution $ a_0 \in (0,2) $.  

Take 
$$
a = \min \left\{a':~ a' \ge \frac{a_0n}{2}, ~ \frac{a'}{4} + \frac{n}{4} ~{\rm is}~ {\rm a}~ {\rm prime}~{\rm number}\right\}.
$$
It is known (see \cite{BHP}) that between $ x $ and $ x + O\left(x^{0.525}\right) $, there is certainly a prime number. Thus, 
the quantity 
$$
p = \frac{a}{4} + \frac{n}{4} 
\eqno{(3)}
$$
is asymptotically equal to $ p_0n $, where 
$$ 
p_0 = \frac{a_0}{8} + \frac{1}{4} \in \left(\frac{1}{4}, \frac{1}{2}\right).
$$

Consider the set 
$$
\Sigma = \{{\bf x} = (x_1, \dots, x_n): ~ \forall~ i ~~ x_i \in \{-1,1\}, ~ x_1 = 1, ~ x_1 + \ldots + x_n = 0\}.
$$ 
Let $ {\cal A} = \left\{a_1, \dots, a_{w}\right\} $ be the set of all possible $2k$-character words over the 
alphabet $ X = \{1, \dots, n\} $, $ w = n^{2k} $. 
Assume that $ a_i = \left\{\nu(i,1), \dots, \nu(i,{2k})\right\} $, $ i = 1, \dots, w $. 
Fix an $ {\bf x} = (x_1, \dots, x_n) \in \Sigma $. Consider 
$$
{\bf x}^{*2k} = \left(x_{\nu(1,1)} \cdot \ldots \cdot x_{\nu(1,2k)}, \ldots, x_{\nu(w,1)} \cdot \ldots \cdot x_{\nu(w,2k)}, 
\sqrt{2k a^{2k-1}} x_1, \ldots, \sqrt{2k a^{2k-1}} x_n\right).
$$
Clearly the number of coordinates in any vector $ {\bf x}^{*2k} $ equals $ w + n $. Put 
$$
\Omega' = \left\{{\bf x}^{*2k}: ~ {\bf x} \in \Sigma\right\}. 
$$
One can readily see that $ \Omega' $ lies in $ {\mathbb R}^{d'} $ with $ d' = w < d $. The point is that for any $ i $,
$ x_i = x_1^{2k-1} \cdot x_i $, since $ x_1 = 1 $. The construction is complete. 

\subsection{Transforming $ \Omega' $ into an $ \Omega \subset S_r^{d-1} $}

First, let us calculate the diameter of $ \Omega' $. For the scalar product of any two vectors $ {\bf x}^{*2k}, {\bf y}^{*2k} \in \Omega' $, 
we have the relation
$$
\left({\bf x}^{*2k},{\bf y}^{*2k}\right) = \sum_{i=1}^w x_{\nu(i,1)} \cdot \ldots \cdot x_{\nu(i,2k)} \cdot 
y_{\nu(i,1)} \cdot \ldots \cdot y_{\nu(i,2k)} + 2k a^{2k-1} ({\bf x},{\bf y}) = 
$$
$$
= \sum_{i_1=1}^n \ldots \sum_{i_{2k}=1}^n x_{i_1} \cdot \ldots \cdot x_{i_{2k}} \cdot y_{i_1} \cdot \ldots \cdot y_{i_{2k}} + 2k a^{2k-1} ({\bf x},{\bf y}) = 
\left(\sum_{i_1=1}^n x_{i_1}y_{i_1}\right) \cdot \ldots \cdot \left(\sum_{i_{2k}=1}^n x_{i_{2k}}y_{i_{2k}}\right) + 2k a^{2k-1} ({\bf x},{\bf y}) = 
$$
$$
= ({\bf x},{\bf y})^{2k} + 2k a^{2k-1} ({\bf x},{\bf y}). 
$$

Obviously the minimum of the form $ \left({\bf x}^{*2k},{\bf y}^{*2k}\right) $ is attained on those and only those pairs of vectors 
$ {\bf x}, {\bf y} \in \Sigma $ whose scalar product equals $ -a $. Such pairs of vectors do really exist for large enough values of $ d $. 
Indeed, by Construction, $ a \sim \frac{a_0n}{2} $ and $ a_0 < 2 $. So for large $ n $, $ -a > -n $. Moreover, $ a \equiv 0 \pmod 4 $, and it is 
easy to see that for every two vectors $ {\bf x}, {\bf y} \in \Sigma $, one necessarily has $ ({\bf x},{\bf y}) \in (-n,n] $ and 
$ ({\bf x},{\bf y}) \equiv 0 \pmod 4 $. 

Thus, we get
$$
\diam^2 \Omega' = 2\left({\bf x}^{*2k},{\bf x}^{*2k}\right) - 2 (a^{2k}-2k a^{2k}) = 2n^{2k} + 4k a^{2k-1} n + (4k-2) a^{2k}. 
$$

At the same time, $ \Omega' $ lies of course on the sphere $ S_{\rho}^{d'-1} $, where 
$$
\rho^2 = \left({\bf x}^{*2k},{\bf x}^{*2k}\right) = n^{2k} + 2k a^{2k-1} n.
$$

Compressing $ \Omega' $ so that a new set $ \Omega'' $ has diameter 1, we see that $ \Omega'' \subset S_{r'}^{d'-1} $ with
$$
(r')^2 = \frac{n^{2k} + 2k a^{2k-1} n}{2n^{2k} + 4k a^{2k-1} n + (4k-2) a^{2k}}.
$$
Since, again by Construction, $ a \ge \frac{a_0n}{2} $ and $ u $ defined in (2) is monotone decreasing, we get the inequalities
$$
(r')^2 = \frac{n^{2k} + 2k a^{2k-1} n}{2n^{2k} + 4k a^{2k-1} n + (4k-2) a^{2k}} \le 
\frac{n^{2k} + 2k \left(\frac{a_0n}{2}\right)^{2k-1} n}{2n^{2k} + 4k \left(\frac{a_0n}{2}\right)^{2k-1} n + 
(4k-2) \left(\frac{a_0n}{2}\right)^{2k}} = u(a_0) = r^2.
$$

Now, it remains to interpret $ S_{r'}^{d'-1} $ as an intersection of the sphere $ S_r^{d-1} $ and a plane of dimension $ d' $. This can be done, 
since $ d' < d $ and $ r' \le r $. Let $ \Omega $ be an image of $ \Omega'' $ under such interpretation. 

\subsection{Lower bound for $ f(\Omega) $}

It is clear that 
$$ 
f_r(d) \ge f(\Omega) = f(\Omega'') = f(\Omega'). 
$$
So it remains to show that $ f(\Omega') \ge (c+\delta)^{\sqrt[2k]{d}} $ for appropriate $ c $ and $ \delta $. First, assume that 
$$ 
f(\Omega') < f = \frac{|\Omega'|}{\sum\limits_{i=0}^{p-1} C_n^i}. 
$$
Then, $ \Omega' $ can be represented as 
$$
\Omega' = \Omega'_1 \cup \ldots \cup \Omega'_f, ~~~ \forall~i ~~ \diam \Omega'_i < \diam \Omega'.
\eqno{(4)}
$$
Obviously the correspondence $ {\bf x} \to {\bf x}^{*2k} $ is a bijection between $ \Sigma $ and $ \Omega' $. So 
$$ 
|\Omega'| = |\Sigma| = 
C_{n-1}^{\frac{n}{2}-1} 
$$ 
and, moreover, the partition (4) induces a partition 
$$
\Sigma = \Sigma_1 \cup \ldots \cup \Sigma_f.
$$
By the choice of $ f $ and by pigeon-hole principle, there exists a part $ \Sigma_i $ with $ |\Sigma_i| > \sum\limits_{i=0}^{p-1} C_n^i $. 
In the next subsection, we shall prove the following lemma. 

\vskip+0.2cm

\noindent {\bf Lemma.} {\it If $ Q \subset \Sigma $ is such that $ |Q| > \sum\limits_{i=0}^{p-1} C_n^i $, then there exist $ {\bf x}, {\bf y} 
\in Q $ with $ ({\bf x},{\bf y}) = -a $.} 

\vskip+0.2cm

Of course Lemma is applied only for $ n \gg 1 $. By this lemma, $ \Sigma_i $ contains two different vectors 
$ {\bf x}, {\bf y} $ with scalar product equal to $ -a $. It means (cf. Subsection 3.2) that 
$$ 
\diam \Omega'_i = |{\bf x}^{*2k}-{\bf y}^{*2k}| = 
\diam \Omega', 
$$ 
which contradicts the properties of partition (4). 

Thus, we have shown that 
$$ 
f(\Omega') \ge \frac{|\Omega'|}{\sum\limits_{i=0}^{p-1} C_n^i} = \frac{C_{n-1}^{\frac{n}{2}-1}}{\sum\limits_{i=0}^{p-1} C_n^i}. 
$$ 

Let us recall that $ p \sim p_0n $, $ p_0 \in \left(\frac{1}{4},\frac{1}{2}\right) $. In this case, standard analytical tools like Stirling's
formula and Chernoff's inequality entail the following asymptotic relations:
$$ 
C_{n-1}^{\frac{n}{2}-1} = (2+o(1))^n, ~~ \sum\limits_{i=0}^{p-1} C_n^i = (c'+o(1))^n, ~~ c' < 2.
$$
Taking the necessary ratio, we get 
$$
f(\Omega') \ge (c+\delta')^n, ~~ c > 1, ~~ \delta' = o(1). 
$$
Finally, by inequalities (1) we have
$$
f_r(d) \ge f(\Omega') \ge (c+\delta')^n \ge (c+\delta')^{\sqrt[2k]{d}-5} = (c+\delta)^{\sqrt[2k]{d}}, ~~ \delta = o(1).
$$
Theorem 2 is proved.   

\subsection{Proof of Lemma}

We shall use the now classical linear algebra method in combinatorics (see \cite{ABS}, \cite{BF}, \cite{Rai5}).

\noindent Consider an arbitrary $ Q = \{{\bf x}_1, \ldots, {\bf x}_s\} \subset \Sigma $ such that for every two different $ i,j $, one has 
$ ({\bf x}_i,{\bf x}_j) \neq -a $. We need to show that $ s \le \sum\limits_{i=0}^{p-1} C_n^i $. 

To each vector $ {\bf x} \in \Sigma $ we assign a polynomial $ P_{{\bf x}} \in {\mathbb Z}/p {\mathbb Z}[y_1, \dots, y_n] $. Namely, 
$$
P_{{\bf x}}({\bf y}) = \prod_{i=0, i \not \equiv -a \pmod p}^{p-1} (i - ({\bf x},{\bf y})), ~~~ {\bf y} = (y_1, \dots, y_n). 
$$
Here the product is taken over all the smallest non-negative residues $ i $ modulo $ p $ except for $ i \equiv -a \pmod p $. So the 
degree of any $ P_{{\bf x}} $ does not exceed $ p-1 $. The most important property of such polynomials is as follows. 

\vskip+0.2cm

\noindent {\bf Property.} {\it For every $ {\bf x}, {\bf y} \in \Sigma $, the congruence $ ({\bf x},{\bf y}) \equiv -a \pmod p $ is 
equivalent to the congruence $ P_{{\bf x}}({\bf y}) \not \equiv 0 \pmod p $.}

\vskip+0.2cm

Property is evident, and we shall just use it. However, before doing so, we make a transformation of any $ P_{{\bf x}} $ into some 
$ P_{{\bf x}}' $. More precisely, we represent every polynomial $ P_{{\bf x}} $ as a linear combination of monomials. Of course 
each monomial has the form 
$$
y_{i_1}^{\alpha_1} \cdot \ldots \cdot y_{i_q}^{\alpha_q}, ~~~ q \le p-1. 
$$
If $ \alpha_{\nu} $ is even, we remove $ y_{i_{\nu}} $ from the monomial. If it is odd, we replace it by 1. Clearly the new monomial 
is just a product of some variables. However, the new polynomials $ P_{{\bf x}}' $, $ {\bf x} \in \Sigma $, still are subject to Property. 
The point is that in Property, only variables whose values are $ \pm 1 $ are considered. 

By the just-given construction, 
$$
{\rm dim} ~\left({\rm linear} ~ {\rm span} \, \left\{P_{{\bf x}}'\right\}_{{\bf x} \in \Sigma}\right) \le \sum\limits_{i=0}^{p-1} C_n^i. 
$$
If we succeed now in showing that the vectors $ {\bf x}_1, \ldots, {\bf x}_s $ from the set $ Q $ correspond to the set 
$ P_{{\bf x}_1}', \dots, P_{{\bf x}_s}' $ of linearly independent polynomials (over the field $ {\mathbb Z}/p {\mathbb Z} $), then 
Lemma is proved. 

So let us assume that 
$$
c_1 P_{{\bf x}_1}'({\bf y}) + \ldots +  c_s P_{{\bf x}_s}'({\bf y}) \equiv 0 \pmod p ~~ \forall~ {\bf y} \in \Sigma. 
\eqno{(5)}
$$
Take $ {\bf y} = {\bf x}_i $ with an arbitrary $ i $. On the one hand, $ ({\bf x}_i,{\bf x}_i) = n $. By Construction, $ n-4p = -a $
(see (3)). Therefore, $ ({\bf x}_i,{\bf x}_i) \equiv -a \pmod p $, i.e., by Property, $ P_{{\bf x}_i}'({\bf x}_i) \not \equiv 0 \pmod p $.
On the other hand, if $ j \neq i $, then $ ({\bf x}_i, {\bf x}_j) < n $ and $ ({\bf x}_i, {\bf x}_j) \neq -a $. Moreover, since 
$$ 
a \equiv 0 \pmod 4, ~~ n \equiv 0 \pmod 4, ~~ ({\bf x},{\bf y}) \equiv 0 \pmod 4 ~~ \forall~ {\bf x},{\bf y} \in \Sigma,
$$
we see that 
$$
({\bf x}_i, {\bf x}_j) \not \in \{n-p, n-2p, n-3p, n-5p, n-6p, n-7p\}.
$$
Finally, $ -a < 0 $, and so $ n-8p < -n $, which means that $ ({\bf x}_i, {\bf x}_j) \not \equiv a \pmod p $ and that, by Property, 
$ P_{{\bf x}_j}'({\bf x}_i) \equiv 0 \pmod p $.  

By relation (5), we get $ c_i \equiv 0 \pmod p $ (here the primality of $ p $ is essential). Since $ i $ was arbitrary, we obtain 
the linear independence of our polynomials, and Lemma is proved.

\section{Proof of Theorem 3}

Take $ r = \frac{1}{2} + \varphi $, where $ \varphi(d) = \frac{6 \ln\ln d}{\ln d} $. We shall prove that for large enough $ d $, 
$ f_r(d) > d+2 $. Then, the whole assertion of Theorem 3 will follow, since for any $ r' > r $, $ S_{r'}^{d} \supset S_r^{d-1} $.

Of course, further exposition will be very close to the one in Section 3. However, there will be some technical subtleties. First of 
all, we rewrite Subsection 3.1. 

Take $ k = \lceil 1/\varphi \rceil $. Then, for large values of $ d $, 
$$
\frac{2k+1}{8k} = \frac{1}{4} + \frac{1}{8k} \le \frac{1}{4} + \frac{\varphi}{8} < \frac{1}{4} + \varphi + \varphi^2 = r^2. 
$$
In other similar inequalities below, we shall not write ``for large values of $ d $'' anymore; we shall just assume that $ d $ is big enough. 
As in Subsection 3.1, put 
$$
n = \max\left\{m: ~ m \equiv 0 \pmod 4, ~~ m^{2k} < d\right\}. 
$$
Now, inequality (1) means that 
$$
n \sim d^{\frac{\varphi}{2}} = e^{(\ln d) \cdot \frac{3 \ln\ln d}{\ln d}} = \ln^3 d.
$$

Let $ a_0 = 2 - \frac{\varphi}{2} $. We want to show that $ u(a_0) < r^2 $, provided $ u $ is the same function as in (2). A standard 
computation is below. First,
$$
u(a_0) = \frac{1+2k\left(\frac{a_0}{2}\right)^{2k-1}}{2+4k \left(\frac{a_0}{2}\right)^{2k-1} + 
(4k-2) \left(\frac{a_0}{2}\right)^{2k}} < \frac{2k (a_0/2)^{2k-1}}{(8k-2) (a_0/2)^{2k}} \cdot 
\frac{1+\frac{1}{2k (a_0/2)^{2k-1}}}{1+\frac{2}{(8k-2) (a_0/2)^{2k}}}.
$$
Further, 
$$
\frac{2k}{8k-2} = \frac{1}{4} \cdot \frac{1}{1-\frac{1}{4k}}. 
$$
Since $ \frac{1}{1-x} \le 1+2x $, we have 
$$ 
\frac{2k}{8k-2} \le \frac{1}{4} \cdot \left(1+\frac{\varphi}{2}\right).
$$ 
Now, 
$$
\frac{(a_0/2)^{2k-1}}{(a_0/2)^{2k}} = \frac{a_0}{2} = 1 - \frac{\varphi}{4}.
$$
Finally, 
$$
\frac{1+\frac{1}{2k (a_0/2)^{2k-1}}}{1+\frac{2}{(8k-2) (a_0/2)^{2k}}} \le 
\left(1+\frac{1}{2k (a_0/2)^{2k-1}}\right) \cdot \left(1 - \frac{1}{(8k-2) (a_0/2)^{2k}}\right).
\eqno{(6)}
$$ 
Clearly 
$$
\left(\frac{a_0}{2}\right)^{2k} \sim \left(\frac{a_0}{2}\right)^{2k-1} \sim \left(1-\frac{\varphi}{4}\right)^{\frac{2}{\varphi}} \sim 
e^{-1/2}. 
$$
Consequently, the right-hand side of (6) can be bounded from above by 
$$
\left(1+\frac{1.1 \varphi}{2 e^{-1/2}}\right) \cdot \left(1-\frac{0.9 \varphi}{8 e^{-1/2}}\right) \le 
1+\frac{3.6 \varphi}{8 e^{-1/2}}.
$$
Thus, 
$$
u(a_0) \le \frac{1}{4} \cdot \left(1+\frac{\varphi}{2}\right) \cdot \left(1 - \frac{\varphi}{4}\right) \cdot \left(1+\frac{3.6 \varphi}{8 e^{-1/2}}\right) 
\le \frac{1}{4} \cdot \left(1 + \frac{\varphi}{4}\right) \cdot \left(1+\frac{3.6 \varphi}{8 e^{-1/2}}\right) \le \frac{1}{4} + \frac{\varphi}{4} 
< \frac{1}{4} + \varphi + \varphi^2 = r^2.
$$

Again, as in Construction, we take 
$$
a = \min \left\{a':~ a' \ge \frac{a_0n}{2}, ~ \frac{a'}{4} + \frac{n}{4} ~{\rm is}~ {\rm a}~ {\rm prime}~{\rm number}\right\}.
$$
By the already cited results from \cite{BHP}, $ a < n $ and 
$$
p = \frac{a}{4} + \frac{n}{4} \le \frac{n}{2} - \frac{\varphi n}{16} + O\left(n^{0.525}\right). 
$$
On the one hand, 
$$
\varphi n = \Theta\left((\ln\ln d) (\ln d)^2\right).
$$
On the other hand, 
$$
n^{0.525} = (\ln d)^{3 \cdot 0.525} = (\ln d)^{1.575} = o\left((\ln\ln d) (\ln d)^2\right). 
$$
Hence, 
$$
p \le \frac{n}{2} - \frac{\varphi n}{20}.
$$

The sets $ \Sigma $ and $ \Omega' $ are just the same as in Subsection 3.1. The only changement in Subsection 3.2, is in replacing the 
last equality by the inequality, in the estimate of $ (r')^2 $ by $ r^2 $. 

Reproducing Subsection 3.3 word by word, we obtain the bound 
$$
f_r(d) \ge \frac{C_{n-1}^{\frac{n}{2}-1}}{\sum\limits_{i=0}^{p-1} C_n^i} = 
\frac{\frac{1}{2} C_{n}^{\frac{n}{2}}}{\sum\limits_{i=0}^{p-1} C_n^i}.
$$

It remains to establish the estimate
$$
\frac{\frac{1}{2} C_{n}^{\frac{n}{2}}}{\sum\limits_{i=0}^{p-1} C_n^i} > d+2.
$$
Clearly 
$$
\sum\limits_{i=0}^{p-1} C_n^i < n C_n^p \le n C_n^{\frac{n}{2}-x}, ~~~ x = \left[\frac{\varphi n}{20}\right]. 
$$
We have
$$
\frac{C_{n}^{\frac{n}{2}}}{C_n^{\frac{n}{2}-x}} = \frac{\left(\frac{n}{2}+1\right) \cdot \ldots \cdot \left(\frac{n}{2}+x\right)}{
\left(\frac{n}{2}\right) \cdot \left(\frac{n}{2}-1\right) \cdot \ldots \cdot \left(\frac{n}{2}-x+1\right)} = 
\frac{\left(1+\frac{2}{n}\right) \cdot \ldots \cdot \left(1+\frac{2x}{n}\right)}{
\left(1-\frac{2}{n}\right) \cdot \ldots \cdot \left(1-\frac{2(x-1)}{n}\right)} = 
$$
$$
= e^{\frac{x(x+1)}{2} - \frac{x(x-1)}{2} + O\left(\frac{x^3}{n^2}\right)} = e^{\frac{2x^2}{n} + O\left(\frac{x^3}{n^2}\right)} = 
e^{(1+o(1))\frac{\varphi^2 n}{200}} = e^{\Theta\left(\frac{(\ln\ln d)^2}{(\ln d)^2} \cdot (\ln d)^3\right)} = 
e^{\Theta\left((\ln\ln d)^2 \cdot (\ln d)\right)}.
$$
Thus, 
$$
\frac{\frac{1}{2} C_{n}^{\frac{n}{2}}}{\sum\limits_{i=0}^{p-1} C_n^i} = d^{\Theta\left((\ln\ln d)^2\right)} > d+2, 
$$
and the proof is complete. 

\section{Proof of Theorem 4}

Take $ S = S_r^{d-1} $ with an arbitrary $ r $. Inscribe into $ S $ a $ d $-dimensional regular simplex $ \Delta $. Consider 
its faces $ \Delta_1, \dots, \Delta_{d+1} $. If $ O $ is the origin, then denote by $ S_i $ the set 
$$
{\rm conv}~ \{\Delta_i,O\} \cap S.
$$
Of course 
$$
S = S_1 \cup \ldots \cup S_{d+1}.
$$
Moreover, it is well-known that (see \cite{BMS})
$$
\diam S_i = 2r\left(1-\Theta\left(\frac{1}{d}\right)\right).
$$
If $ r = r(d) \le \frac{1}{2} + \frac{c}{d} $ with an appropriate $ c > 0 $, then we get $ \diam S_i < 1 $, and we are done. 

\medskip

In principle, it's possible to obtain rather good upper bounds for $ f_r(d) $ in general case. The simplest way for doing that is to use old results 
of C.A. Rogers (see \cite{Rog}): {\it any sphere of radius $ r > \frac{1}{2} $ in $ {\mathbb R}^d $ can be covered by $ (2r+o(1))^d $
spherical caps of diameter 1.} This result gives already the estimate $ f_r(d) \le (2r+o(1))^d $ when $ r $ is fixed. More subtle 
estimates can be discovered by using an approach which appeared in the paper \cite{BL}. For example, Rogers tells us that for 
$ r > \sqrt{\frac{3}{8}} $, we can only show that $ f_r(d) \le (c+o(1))^d $ with $ c > \sqrt{\frac{3}{2}} $. However, J. Bourgain and 
J. Lindenstrauss, the authors of \cite{BL}, provide us with a universal bound $ f_r(d) \le \left(\sqrt{\frac{3}{2}}+o(1)\right)^d $. Their 
ideas may be carefully applied even in some cases when $ r \to \frac{1}{2} $, but, in this paper, we do not dwell on that kind of results. 

\section{Improving Construction from Subsection 3.1 is hard}

Let us discuss some key properties of Construction from \S 3.1 and a possibility of improving it. Indeed, one of the most important steps in 
Construction was in assigning, to each $ n $-dimensional vector $ {\bf x} \in \Sigma $, a $ d' $-dimensional vector $ {\bf x}^{*2k} $, so that 
eventually we got a set $ \Omega' \subset {\mathbb R}^{d'} $. The correspondence between vectors $ {\bf x} $ from $ \Sigma $ and their images 
$ {\bf x}^{*2k} $ was organized in such a way that 
$$
\left({\bf x}^{*2k},{\bf y}^{*2k}\right) = ({\bf x},{\bf y})^{2k} + 2k a^{2k-1} ({\bf x},{\bf y}).
$$
In other words, the scalar product of any two elements of the set $ \Omega' $ was a polynomial depending on the scalar product of the 
preimages of those elements. Denote this polynomial by $ h_a^* $. Thus, $ h_a^*(t) = t^{2k} + 2ka^{2k-1} t $. 

In order to prove Theorems 1 -- 3, we essentially used the following properties of $ h_a^* $ (see \S 3.2): first of all, $ h_a^* $, considered as  
function on the interval $ [-n,n] $, attains its minimum at point $ -a \in (-n,0) $ and so 
$$ 
\left.\frac{{\rm d} h_a^*}{{\rm d} t}\right|_{t=-a} = 0; 
$$
second, the quantity
$$
\frac{h_a^*(n)}{2h_a^*(n)-2h_a^*(-a)} = \frac{n^{2k} + 2k a^{2k-1} n}{2n^{2k} + 4k a^{2k-1} n + (4k-2) a^{2k}} 
$$
is a close approximation to the radius of a sphere, on which a counterexample to Borsuk's conjecture lies. Therefore, the results of Theorems 
1 -- 3 would be improved, provided we could replace $ h_a^* $ by another polynomial $ h_a $ having the same first property and a smaller value of 
the expression 
$$
\frac{h_a(n)}{2h_a(n)-2h_a(-a)}. 
$$
Moreover, we may also vary the value $ a \in (0,n) $ by taking $ a \sim a'n $ with $ a' \in (0,1) $ as close to 1 as necessary.

Let $ \Delta_a^m $ be the set of all polynomials $ h_a(t) = u_m t^m + \ldots + u_1 t + u_0 $ of degree $ m $ with non-negative coefficients 
and having the property that $ h_a $, considered as function on the interval $ [-n,n] $, attains its minimum at point $ -a \in [-n,0) $ and  
$$ 
\left.\frac{{\rm d} h_a}{{\rm d} t}\right|_{t=-a} = 0
$$
(the last condition is non-trivial only for $ a = n $). 

\vskip+0.2cm

\noindent {\bf Proposition.} {\it If $ m $ is even, then 
$$
\min_{a \in (0,n]} \min_{h_a \in \Delta_a^m} \frac{h_a(n)}{2h_a(n)-2h_a(-a)} = 
\min_{h \in \Delta_n^m} \frac{h(n)}{2h(n)-2h(-n)}
= \frac{h^*(n)}{2h^*(n)-2h^*(-n)},
$$
where $ h^*(t) = t^m + m n^{m-1} t $. If $ m $ is odd, then 
$$
\min_{h \in \Delta_n^m} \frac{h(n)}{2h(n)-2h(-n)} \ge 
\min_{h \in \Delta_n^{m-1}} \frac{h(n)}{2h(n)-2h(-n)}.
$$}

\vskip+0.2cm

Before proving Proposition, let us briefly comment on it. Actually, Proposition tells us that, in \S 3.1 and \S 3.2, everything was done 
in an optimum way: it's better to take polynomials of even degree, and, among them, $ h_a^* $ is asymptotically best possible
(when $ a \to n $, $ h_a^* \to t^{2k} + 2k n^{2k-1} t $). It is worth noting that, 
in principle, for any $ h \in \Delta_a^m $, one can transform $ \Sigma $ into such an $ \Omega' $ that the scalar product of any two vectors 
from $ \Omega' $ is equal to the value of $ h $ at the scalar product of the preimages of those two vectors (here it is important to assume 
that any $ h \in \Delta_a^m $ has non-negative coefficients). However, according to 
Proposition, this fact is already not quite useful for our purposes.   
   
\paragraph{Proof of Proposition.} First, we note that the quantity $ h_a(-a) $ should be non-positive in order to minimize the ratio 
$$
\frac{h_a(n)}{2h_a(n)-2h_a(-a)} = \frac{1}{2-\frac{2h_a(-a)}{h_a(n)}}, 
\eqno{(7)}
$$
since $ h_a(n) > 0 $ for any $ h_a $. Polynomials $ h_a \in \Delta_a^m $ such that $ h_a(-a) \le 0 $ do really exist, so we may assume that 
$ h_a(-a) \le 0 $. Under this assumption, our minimization is equivalent to the maximization of the expression $ \frac{|h_a(-a)|}{h_a(n)} $
over the set $ \tilde{\Delta}_a^m \subset \Delta_a^m $ containing those and only those polynomials $ h_a $, for which $ h_a(-a) \le 0 $. 

Furthermore, we may suppose that $ u_0 = 0 $, since, for $ u_0 > 0 $, ratio (7) is definitely greater. 

Now, let us prove that 
$$
\min_{a \in (0,n]} \min_{h_a \in \Delta_a^m} \frac{h_a(n)}{2h_a(n)-2h_a(-a)} = 
\min_{h \in \Delta_n^m} \frac{h(n)}{2h(n)-2h(-n)}
\eqno{(8)} 
$$ 
or, which is the same, that 
$$
\max_{a \in (0,n]} \max_{h_a \in \tilde{\Delta}_a^m} \frac{|h_a(-a)|}{h_a(n)} = 
\max_{h \in \tilde{\Delta}_n^m} \frac{|h(-n)|}{h(n)}.
\eqno{(9)} 
$$ 
Indeed, all the coefficients of $ h_a $ are non-negative, and so $ h_a(a) \le h_a(n) $. Therefore, for any $ a \in (0,n] $, we have
$$
\max_{h_a \in \tilde{\Delta}_a^m} \frac{|h_a(-a)|}{h_a(n)} \le \max_{h_a \in \tilde{\Delta}_a^m} \frac{|h_a(-a)|}{h_a(a)} = 
\max_{h \in \tilde{\Delta}_n^m} \frac{|h(-n)|}{h(n)},
$$
which completes the proof of (8) and (9). 

It remains to show that 
$$
\max_{h \in \tilde{\Delta}_n^m} \frac{|h(-n)|}{h(n)} \le \frac{|h^*(-n)|}{h^*(n)} = \frac{m-1}{m+1}.
$$
Take an arbitrary polynomial $ h \in \tilde{\Delta}_n^m $. We represent it in the form 
$$
h = h_1+h_2 = \sum_i c_{\alpha_i} t^{\alpha_i}+\sum_i c_{\beta_i} t^{\beta_i}, ~~~ c_{\alpha_i} > 0, ~~ c_{\beta_i} > 0,
$$
where, in $ h_1 $, only odd degrees of $ t $ are taken, and, in $ h_2 $, only even degrees are present. In this notation, we have 
$$
\frac{|h(-n)|}{h(n)} = \frac{\left|\sum\limits_i c_{\alpha_i} n^{\alpha_i}-\sum\limits_i c_{\beta_i} n^{\beta_i}\right|}{\sum\limits_i c_{\alpha_i} 
n^{\alpha_i}+\sum\limits_i c_{\beta_i} n^{\beta_i}}.
$$  
Hence, we are led to check the inequality 
$$
\sum\limits_i c_{\alpha_i} n^{\alpha_i} \le m \sum\limits_i c_{\beta_i} n^{\beta_i}.
\eqno{(10)}
$$
However, we know that $ h'(-n) = 0 $, i.e., 
$$
\sum\limits_i \alpha_i c_{\alpha_i} n^{\alpha_i-1} = \sum\limits_i \beta_i c_{\beta_i} n^{\beta_i-1}.
$$
Then, we get the following series of inequalities:
$$
\sum\limits_i c_{\alpha_i} n^{\alpha_i} \le n \sum\limits_i \alpha_i c_{\alpha_i} n^{\alpha_i-1} = n \sum\limits_i \beta_i c_{\beta_i} n^{\beta_i-1}
\le m \sum\limits_i c_{\beta_i} n^{\beta_i}.
$$ 
Thus, (10) is true, and the proof of the first part of Proposition is complete. 

\vskip+0.2cm

The second part would follow from the inequality 
$$
\max_{h \in \tilde{\Delta}_n^m} \frac{|h(-n)|}{h(n)} \le \frac{m-2}{m},
$$
which is tantamount to 
$$
\sum\limits_i c_{\alpha_i} n^{\alpha_i} \le (m-1) \sum\limits_i c_{\beta_i} n^{\beta_i}.
\eqno{(11)}
$$
Eventually, (11) is provided by the following series of estimates:
$$
\sum\limits_i c_{\alpha_i} n^{\alpha_i} \le n \sum\limits_i \alpha_i c_{\alpha_i} n^{\alpha_i-1} = n \sum\limits_i \beta_i c_{\beta_i} n^{\beta_i-1}
\le (m-1) \sum\limits_i c_{\beta_i} n^{\beta_i}.
$$ 
Here we get the factor $ m-1 $ instead of $ m $, since $ m $ is odd and so $ \beta_i \le m-1 $ for any $ i $. 

Proposition is proved.

\section{Improving Construction from Subsection 3.1 is still possible}

The only place in Construction, which may me further improved, is in the form of the set $ \Sigma $. In principle, we could take 
an arbitrary subset $ \Sigma $ of the integer lattice $ {\mathbb Z}^n $. In order to be able to apply the whole machinery we used to prove 
our results we must assume that for any $ {\bf x} \in \Sigma $, the value of $ ({\bf x},{\bf x}) $ is the same. For example, 
$$
\Sigma = \{{\bf x} = (x_1, \dots, x_n): ~ \forall ~ i ~~ x_i \in \{-1,0,1\}, ~ x_1^2 + \ldots + x_n^2 = k\}
$$
with some $ k $ could be of help. In the papers \cite{Rai6}, \cite{Rai7} and in the book \cite{Rai5}, an analogous approach was developed 
in order to get optimal bounds for the chromatic numbers of spaces and for the standard Borsuk number. However, the difference between the 
results obtained from (-1,1)-constructions and the more general ones turned out to appear only in the expressions of the form $ c+\delta $ (cf. 
Theorem 2). Although sometimes this is important as well, we do not think it is of interest to thoroughly investigate here the 
corresponding bounds.   

Another possible refinement is in calculating more carefully the dimension of the set $ \Omega' $. Actually, it is not $ n^{2k} $; it does not 
exceed $ C_n^{2k} \le \frac{n^{2k}}{(2k)!} $. For fixed values of $ r > 1/2 $, this is very important, but for sequences $ r_d \sim 1/2 $, 
this fact does not give anything.

%\end{large}

\end{document}